\documentclass[a4paper,10pt]{article}

\def\endproof{\mbox{\ \rule{.1in}{.1in}}}
\usepackage{amsfonts}
\usepackage{graphicx}
\usepackage{dsfont}
\newtheorem{claim}{Claim}
\newtheorem{lemma}{Lemma}

\title{Balanced Convex Partitions of Measures in $\mathds{R}^d$}
\author{Pablo Sober\'on \\ pablo@math.unam.mx \and ------------------------------------------------------------ \\ Department of Mathematics \\ University College London \\ Gower Sreet, London WC1E 6BT \\ United Kingdom
\\ ------------------------------------------------------------}
\date{}
\begin{document}

\maketitle

\begin{abstract}
We will prove the following generalization of the ham sandwich Theorem, conjectured by Imre B\'ar\'any.  Given a positive integer $k$ and $d$ nice measures $\mu_1, \mu_2, \ldots, \mu_d$ in $\mathds{R}^d$ such that $\mu_i (\mathds{R}^d) = k$ for all $i$, there is a partition of $\mathds{R}^d$ in $k$ interior-disjoint convex parts $C_1, C_2, \ldots, C_k$ such that $\mu_i (C_j) = 1$ for all $i,j$.  If $k=2$ this gives the ham sandwich Theorem.  This result was independently proved by R.N. Karasev.
\end{abstract}

\section{Introduction}

The ham sandwich Theorem is a very well known result in both measure theory and discrete geometry (see \cite{Bey04} and the references therein).  It says the following.
\\

\noindent \textbf{Theorem A (Ham Sandwich)}
Given $d$ finite measures $\mu_1, \mu_2, \ldots, \mu_d$ in $\mathds{R}^d$ that vanish on every hyperplane, there is a hyperplane $H$ that simultaneously divides equally all the measures.
\\

In the planar case, Toshinori Sakai proved that if one wants to partition the measures in more parts, one can find a partition of the plane in convex sets with this property \cite{Sak02}.  The same result was proven by Bespamyatnikh, Kirkpatrick and Snoeyink for two discrete sets of points.  Namely,
\\

\noindent \textbf{Theorem B (Bespamyatnikh, Kirkpatrick and Snoeyink, 2000 \cite{Bes00};Sakai, 2002 \cite{Sak02})}
Given a positive integer $k$ and two nice measures in the plane $\mu_1$ and $\mu_2$ such that $\mu_1 (\mathds{R}^2) = \mu_2 (\mathds{R}^2) = k$, there is a convex partition of the plane in sets $C_1, C_2, \ldots, C_k$ such that $\mu_i (C_j) = 1$ for all $i, j$.
\\

A mixed version has also been proven, with convex sets partitioning simultaneously a measure and a discrete set of points \cite{Gun10}.  In Sakai's theorem, a nice measure refers to a measure $\mu$ absolutely continuous with respect to the Lebesgue measure such that there is a bounded domain $B$ with $\mu(B) = \mu(\mathds{R}^2)$.  We will define a nice measure in $\mathds{R}^d$ as a finite measure $\mu$ absolutely continuous with respect to the Lebesgue measure such that there is a bounded convex set $K$ with $\mu(K)=\mu(\mathds{R}^d)$ and $\mu$ is nonvanishing in open sets of $K$.  Notice that if $\mu$ is a nice measure and $K_1$ is a convex body with positive measure, then $\mu |_{K_1}$ is also a nice measure.  Also, since the set of nice measures is dense we are not losing generality.  A convex partition of $\mathds{R}^d$ $(C_1, C_2, \ldots, C_k)$ is a partition of $\mathds{R}^d$ where all the parts are closed convex sets with pairwise disjoint interiors.

The question of whether Sakai's result could be extended to $\mathds{R}^d$ was asked by Imre B\'ar\'any.  Namely, he conjectured the following theorem.
\\

\noindent \textbf{Theorem 1.}
Given a positive integer $k$ and $d$ nice measures $\mu_1, \mu_2, \ldots, \mu_d$ in $\mathds{R}^d$ such that $\mu_i (\mathds{R}^d) = k$ for all $i$, there is a convex partition of $\mathds{R}^d$ in sets $C_1, C_2, \ldots, C_k$ such that $\mu_i (C_j) = 1$ for all $i,j$.
\\

The aim of this paper is to give a positive answer to B\'ar\'any's conjecture.  This theorem is a generalization of the Ham Sandwich theorem since a convex partition $(C_1, C_2)$ of $\mathds{R}^d$ is always defined by a hyperplane.  It was proved indpendently by R.N. Karasev \cite{Kar10} using different topological methods.  The main tools that will be used in the proof are new results regarding power diagrams and a theorem from equivariant topology.  Power diagrams are a generalization of Voronoi diagrams and will be discussed in the next section.  The theorem from equivariant topology we will use is the following.
\\

\noindent \textbf{Theorem C. (Dold, 1983 \cite{Dold83})}  Let $G$ be a finite group, $|G|>1$, let $X$ be an $n$-connected space with a free action of $G$, and let $Y$ be a (paracompact) topological space of dimension at most $n$ with a free action of $G$.  Then there is no $G$-equivariant map $f: X \longrightarrow Y$.
\\

Recall that the Ham Sandwich theorem is a direct consequence of the Borsuk-Ulam Theorem.  The theorem above, with $X= \mathds{S}^{n+1}$, $Y=\mathds{S}^n$ and $G=\mathds{Z}_2$ gives precisely the Borsuk-Ulam Theorem.  Theorem C gives us the necesary freedom to extend the Ham Sandwich result.  The main idea is to construct a space $X$ that represents the partitions and a space $Y$ that represents how they partition each measure.  If there is no equipartition, then the dimension of $Y$ can be reduced and we obtain a contradiction to Dold's theorem.
\\

Topological methods like this one almost always apear in proofs of theorems related to partitions of measures.  See \cite{Mat03} for more details.
\section{Power Diagrams}

Given $S$ an ordered $k$-tuple of different points $(x_1, x_2, \ldots, x_k)$ in $\mathds{R}^d$ (these points will be called sites), and a weight vector $(w_1, w_2, \ldots, w_k) \in \mathds{R}^k$, we will define the power fuctions $h_i: \mathds{R}^d \longrightarrow \mathds{R}$ by $h_i (x) = d(x,x_i)^2 - w_i$.  The power diagram $C(S,w)$ is a partition of $\mathds{R}^d$ where
$$
x \in C_i \Leftrightarrow h_i (x) \le h_j (x) \ \mbox{for all } j
$$
In other words, $C_i$ is the set of points where $h_i$ is minimal among all the power functions.
\\

If $w=(0,0, \ldots, 0)$, then each point is in the part corresponding to the closest site, which is the Voronoi diagram with sites $S$.
\\

Each $C_i$ is the intersection of all the closed halfspaces $H_{i,j}^{+} = \{ x \ | \ h_i (x) \le h_j (x)\}$ for all $j \neq i$.  Thus, each $C_i$ is a convex polyhedron with at most $k-1$ facets.  Notice that the hyperplane 
$$
H_{i,j} = \{ x \ | \ h_i (x)= h_j (x)  \} = \{ x \ | \ 2( x_i-x_j) \cdot x  = w_j - w_i \}
$$
 is orthogonal to $x_i - x_j$, and its position depends entirely on $w_j-w_i$.  Thus, if $v_0 = (1,1, \ldots, 1)$ is the diagonal vector in $\mathds{R}^k$, we have that $C(S,w) = C(S,w+\alpha v_0)$ for all $S, w, \alpha$.  
\\

In 1998, F. Aurenhammer et al. proved \cite{Aro98} that given a nice probability measure $\mu$, a $k$-set $S$ of $\mathds{R}^d$ and a capacity vector $c = (c_1, c_2, \ldots, c_k) \in \mathds{R}^k$ such that $c_i \ge 0$ for all $i$ and the sum of the $c_i$ is $1$, there is a weight vector $w$ such that for the power diagram $C(S,w)$ we have that $\mu(C_i)=c_i$ for all $i$.  It has also been proven that the vector $w$ is unique up to translations by the diagonal \cite{Aro10}.  The nonvanishing condition on our measure is essentially needed here.
\\

This last result will be our main tool to prove the main theorem.  Since $C(S,w) = C(S, w+\alpha v_0)$ for all $\alpha$, we may choose that the dot product $w \cdot c$ is $0$.
\\

It is also known that given $c$, if one moves the points of $S$ continuously and they remain different \cite{Aro10}, then the weight vector also moves continuously (if a condition such as $c \cdot w = 0$ has been imposed).  We will analyze what happens when the points of $S$ move and some of them converge to the same point.
\\

Let $\mu$ be a nice probability measure in $\mathds{R}^d$ and $c= (c_1, c_2, \ldots, c_k)$ be a capacity vector such that $c_i > 0$ for all $i$.  Given $S$ a set of $k$ different sites in $\mathds{R}^d$ let $w(S) = (w_1, w_2, \ldots, w_k)$ be the weight vector such that the measure of the parts of $C(S,w(S))$ agree with the capacity vector and $w(S) \cdot c = 0$.

\begin{claim} \label{continuidad}
 Let $S$ be a $k$-set of sites in $\mathds{R}^d$ that move continuously through the time interval $[0,1]$.  Suppose the sites remain different in the time interval $[0,1)$.  Suppose that there are two points $x_1, x_2 \in S$ and a point $p_0 \in \mathds{R}^d$ such that $x_i \longrightarrow p_0$ as the time $t \longrightarrow 1$ for $i = 1,2$.  Then, there is a number $w'$ such that $w_i \longrightarrow w'$ as the time $t \longrightarrow 1$ for $i=1,2$.
\end{claim}

In other words, if two sites converge to the same point, their weights converge to the same number.
\\

\noindent \textbf{Proof}.  Suppose there are numbers $w_1 '$ and $w_2 '$ such that $w_i \longrightarrow w_i'$ for $i=1,2$ and $w_1' > w_2'$.  Let $y$ be the point where the hyperplane $H_{1,2} = \{ x \ | \ d(x,x_1)^2 - d(x,x_2)^2 = w_1 - w_2 \}$ intersects the line that goes through $x_1$ and $x_2$.  Let $u=d(y,x_2), v=d(x_1,x_2)$.  Since $w_i \longrightarrow w_i'$ we can suppose that $w_1 > w_2$, so $d(y,x_1) = u+v$.  Thus $w_1 - w_2 = (u+v)^2 - u^2 = v(2u + v)$.  Notice that $(w_1 - w_2) \longrightarrow (w_1' - w_2') > 0$ and $v \longrightarrow 0$.  Thus $u \longrightarrow \infty$.  This means that the distance $d(p_0, C_2) \longrightarrow \infty$.  This contradicts the fact that $\mu (C_2) = c_2 >0$ for $t \in [0,1)$.

Now we need to show that such limit exists.  Using the same argument and the fact that $w \cdot c = 0$, one can see that the weight vector must be bounded as $t$ approaches $1$.  Thus, it suffices to show that if any sequence of values of $w$ as the time approaches $1$ converges, it does to the same limit.  Let $w^{(1)}, w^{(2)}, \ldots$ be such a sequence and $\tilde{w}$ be its limit.  We will replace $S, c, \tilde{w}$ by $S', c', \tilde{w}'$ in the following way:
\begin{itemize}
 \item $S'$ consists of replacing all the points converging to $q$ by a single copy of $q$, for all $q \in \mathds{R}^d$.
 \item $c'$ consists of replacing all the capacities corresponding the points converging to $q$ by a single copy of their sum, for all $q \in \mathds{R}^d$.
 \item $\tilde{w}'$ consists of replacing all the weights corresponding to the points convergint to $q$ by a single copy of that number (we already proved that it must be the same), for all $q \in \mathds{R}^d$. 
\end{itemize}

If we do this, it is clear that the parts of the power diagram of $C(S',w')$ have measures that agree with $c'$.  Also, $\tilde{w}' \cdot c' = 0$.  Thus, there is only one possible value for $\tilde{w}'$, and the same happens for $\tilde{w}$.\endproof

\section{Main proofs}

We will now construct the spaces and functions to use Dold's theorem.  The main idea goes as follows.  Let $p$ be a prime number and let $X$ be the space of ordered $p$-tuples of vectors of $\mathds{R}^d$ such that they are not all the same point.  Notice that $X \cong \mathds{R}^{pd} \backslash Y$ where $Y \cong \mathds{R}^d$.  The space $X$ will be used to represent the partitions of $\mathds{R}^d$ in at least two parts (one for each different point in our $p$-tuple).  The measure we want each part to have will be corresponding with how many times its point is present in the $p$-tuple.  We will use Dold's theorem to see that this can be done simoultaneously for all measures.  The reason why $p$ is required to be prime is only to make sure that  the group actions we will define are free.
\\

\noindent The following lemma will be the core of our proof of the main theorem.

\begin{lemma}
 Let $p$, $d$ be positive integers such that $p$ is prime.  Let $\mu_1, \mu_2, \ldots, \mu_d$ be $d$ nice measures in $\mathds{R}^d$ such that $\mu_i (\mathds{R}^d) = p$ for all $i$.  Then, there is an integer $2 \le r \le p$ and a partition of $\mathds{R}^d$ in $r$ convex parts $C_1, C_2, \ldots, C_r$ such that $\mu_i(C_j) = \mu_{i'} (C_j)$ for all $i,i',j$ and all these measures are positive integers.
\end{lemma}

\noindent \textbf{Proof.}
Given $1 \le i \le d$, let us define a function $f_i : X \longrightarrow \mathds{R}^p$ (associated with $\mu_i$).  Given $x= (x_1, x_2, \ldots, x_p) \in X$, let $S(x) = (s_1, s_2, \ldots, s_t)$ be the $t$-tuple of different points in $x$, with the order in which they appeared in $x$.  Notice that if all the points of $x$ are different, then $S(x) = x$.  For each $1 \le j \le t$, let $\alpha_j$ be the number of times $s_j$ appeared in $x$.  Then, let $(w_1, w_2, \ldots, w_t)$ be the weights needed to partition $\mu_i$ with a power diagram with sites $S(x)$ and capacities $c=(\alpha_1, \alpha_2, \ldots, \alpha_t)$ such that $w \cdot c = 0$.  For all $1 \le j \le p$, let $y_j = w_h$ if $x_j = s_h$.  Now define $f_i (x) = (y_1, y_2, \ldots, y_p)$.  Since $c \cdot w = 0$, we have that $y_1 + y_2 + \ldots + y_p = 0$.  Thus, $f_i(x) \in \mathds{R}^{p-1} \hookrightarrow \mathds{R}^p$.  By Claim \ref{continuidad}, each $f_i$ is continuous.
\\

Now consider $f: X \longrightarrow \mathds{R}^{(p-1)(d-1)} \hookrightarrow \mathds{R}^{p(d-1)}$ defined by $f= (f_1 - f_2, f_1 - f_3, \ldots, f_1 - f_d)$.  We will show that there is an $x \in X$ such that $f(x) = 0$.
\\

Suppose there is no such $x$, so $f: X \longrightarrow \mathds{R}^{(p-1)(d-1)} \backslash \{ 0 \}$.  There is a natural action of $\mathds{Z}_p$ in $X$ given by $\sigma (x_1, x_2, \ldots, x_p) = (x_2, x_3, \ldots, x_p, x_1)$, where $\sigma$ is a generator of $\mathds{Z}_p$.  The same action can be applied in $\mathds{R}^{p-1} = \{ (z_1, z_2, \ldots, z_p) \ | \ z_1 + z_2 + \ldots + z_p = 0 \}$ and thus in $\mathds{R}^{(p-1)(d-1)}$.  Notice that with these actions, $f$ is equivariant.  Since $p$ is prime, the actions in $X$ and in $\mathds{R}^{(p-1)(d-1)} \backslash \{ 0 \}$ are both free.

Let $g(x) = \frac{f(x)}{||f(x)||}$.  We know that $g: X \longrightarrow \mathds{S}^{(p-1)(d-1)-1}$ and is still equivariant.  However, $X$ is $\mathds{R}^{pd}$ with a hole of dimension $d$, so it is at least $[(pd-1) - (d+1)]$-connected.  Since $(pd-1)-(d+1) \ge (p-1)(d-1)-1$, we have a contradictions with Theorem $C$.
\\

Given the point $x$ such that $f(x) = 0$, we have that for $S(x)$ and its associated capacity vector $c$ as in the definition of $f_i$, the weight vector $w$ such that $C(S(x),w)$ partitions $\mu_i$ with capacities $c$ is the same for all $i$.  Thus $C(S(x),w)$ is the partition we wanted.\endproof

\begin{lemma}
Let $a,b$ be positive integers. If Theorem $1$ is true for values $k=a$ and $k=b$, then it is true for $k=ab$.
\end{lemma}

\noindent \textbf{Proof.}
Let $\mu_1, \mu_2, \ldots, \mu_d$ be nice measures such that $\mu_i (\mathds{R}^d) = ab$ for all $i$.  Since Theorem $1$ is true for $k=a$, we can find a convex partition $(C_1, C_2, \ldots, C_a)$ such that $\mu_i (C_j) = b$ for all $i,j$.  For all $1 \le j \le a$, $\mu_i | _{C_j}$ is a nice measure in $\mathds{R}^d$.  Thus, since Theorem $1$ holds for $k=b$, we can find a partition of $\mathds{R}^d$ in convex sets $C_{j,1}, C_{j,2}, \ldots, C_{j,b}$ such that $\mu_i |_{C_j} (C_{j,h}) = 1$ for all $i,h$.  Let $D_{j,h} = C_j \cap C_{j,h}$.  We have that $\mu_i (D_{j,h}) = \mu_i |_{C_j} (C_{j,h}) = 1$.  Since $D_{j,h}$ is convex for all $j,h$, they form the partition of $\mathds{R}^d$ we were looking for.\endproof
\\

We are now ready to prove the main theorem.
\\

\noindent \textbf{Theorem 1.}
Given a positive integer $k$ and $d$ nice measures $\mu_1, \mu_2, \ldots, \mu_d$ in $\mathds{R}^d$ such that $\mu_i (\mathds{R}^d) = k$ for all $i$, there is a convex partition of $\mathds{R}^d$ in sets $C_1, C_2, \ldots, C_k$ such that $\mu_i (C_j) = 1$ for all $i,j$.
\\

\noindent \textbf{Proof.}
We will prove this theorem by strong induction on $k$.  For the basis, if $k=1$ then $C_1 = \mathds{R}^d$.  Suppose that $k \ge 2$ and the theorem is true for all $1 \le k' < k$.
\begin{itemize}
 \item If $k$ is not prime, it can be factorized as $k=ab$ with $a,b < k$.  By lemma $2$, we are done.
 \item If $k$ is prime, by Lemma $1$, one can find a partition of $\mathds{R}^d$ in at least $2$ convex sets $C_1, C_2, \ldots, C_r$ such that each convex set has the same measure in all $d$ measures and these numbers are all positive integers.  Given $1\le j \le r$, we can apply Theorem $1$ in the measures $\mu_i |_{C_j}$ as in the proof of Lemma $2$ to obtain the partition of $\mathds{R}^d$ we were looking for. \endproof
\end{itemize}

\noindent Since the set of nice measures is dense, by usual approximation arguments one can make one of the measures to be a Dirac measure in the origin.  By doing this, we obtain the following Corollary.
\\

\noindent \textbf{Corollary 1.}
Given $\mu_1, \mu_2, \ldots, \mu_{d-1}$ nice measures on $\mathds{S}^{d-1}$ such that $\mu_i (\mathds{S}^{d-1}) = k$ for all $i$, there is a convex cone subdivision (with apices at the origin) $C_1, C_2, \ldots, C_k$ such that $\mu_i (C_j) = 1$ for all $i,j$.
\\

Thus the main theorem also holds for measures in $\mathds{S}^d$.  The equivalence with the version in $\mathds{S}^d$ was conjectured and proven by Imre B\'ar\'any.

\section{Remarks on the proof}
 
 The last inequality used in the proof of Lemma 1 is equivalent to $p \ge 2$. If one goes through Sakai's proof of the planar case, he uses a similar lemma. The lemma says the following: \textit{If there is no line that divides both measures equally and each halfspace has a positive integer measure, then there it is possible to do so with a convex $3$-partition}. Moreover, Sakai proves that if the $3$-partition is necessary, then (as it must be a convex $3$-fan) one of the direction of its lines may be chosen arbitrarily.
In the same way for the high-dimensional version, if the smallest partition possible is with t parts, then conditions may be imposed on $X$ so that they do not reduce its connectedness in more than $t-2$.

\bibliographystyle{amsplain}

\bibliography{ref-part.bib}
\end{document}